\documentclass[review]{elsarticle}
\usepackage{amsfonts,amssymb,amsmath,amsbsy,amsthm,amscd}
\usepackage[english]{babel}
\usepackage{graphicx}

\usepackage{color}
\def\br{{\bf r}}
\UseRawInputEncoding
\def\dfrac#1#2{\displaystyle{#1\over #2}}

\def\bv{{\bf V}}
\def\bV{{\bf V}}

\def\bE{{\bf E}}
\def\bY{{\bf Y}}
\def\bS{{\bf S}}
\def\br{{\bf r}}

\def\bx{{\bf x}}
\def\by{{\bf y}}
\def\bE{{\bf E}}

\title{On globally smooth oscillating solutions of  non-strictly hyperbolic systems} 
\author{Olga S. Rozanova} 
\address{Mathematics and Mechanics Department, Lomonosov Moscow State University, Leninskie Gory,
Moscow, 119991,
Russian Federation} 

\theoremstyle{plain}
\newtheorem{theorem}{Theorem}
\newtheorem{lemma}{Lemma}
\newtheorem{proposition}{Proposition}
\newtheorem{definition}{Definition}
\newtheorem{remark}{Remark}

\numberwithin{equation}{section}
\numberwithin{theorem}{section}
\numberwithin{lemma}{section}
\numberwithin{proposition}{section}
\numberwithin{remark}{section}

\begin{document}
\begin{abstract}
A class of non-strictly hyperbolic systems of quasilinear equations with oscillatory solutions of the Cauchy problem, globally smooth in time in some open neighborhood of the zero stationary state, is found. For such systems, the period of oscillation of solutions does not depend on the initial point of the Lagrangian trajectory. The question of the possibility of constructing these systems in a physical context is also discussed, and non-relativistic and relativistic equations of cold plasma are studied from this point of view.
\end{abstract}
\maketitle


\section{Introduction}
Consider a non-strictly hyperbolic system
\begin{eqnarray}\label{main}
\bY_t+{\mathbb A}(x, {\bY}) \, \bY_x &=& {\bS}(x, {\bY}),\quad {\mathbb A}=Q(x, {\bY})\, {\mathbb I},
\end{eqnarray}
where $\bY(t,x)=(Y_1,\dots, Y_n)^T$, $x\in \mathbb R$, $t>0$, ${\bf S}=(S_1,\dots, S_n)^T$, $i=1,\dots, n,$ $n\in \mathbb N$, $\mathbb I$ is the identity matrix $n\times n$, $t>0$, $x\in \mathbb R$, $Q(x, {\bf 0})=0$, ${\bS}(x, {\bf 0})={\bf 0}$,
subject to the initial data
\begin{eqnarray}\label{mainID}
\bY\Big|_{t=0} = \bY^0(x)\in C^2({\mathbb R}).
\end{eqnarray}
It is assumed that the functions $Q(x, {\bY})$ and ${\bS}(x, {\bY})$ are $C^1$-- smooth functions of all their arguments.

The system \eqref{main} is an intermediate mathematical object between quasilinear hyperbolic equations of general form, where $\mathbb A$ is an arbitrary $n\times n$ matrix with $n$ real distinct eigenvalues, and systems of nonlinear ordinary differential equations. Indeed, the dynamics of the solution can be completely described by the behavior along a single Lagrangian characteristic $x=x(t)$, which is governed by a system of $n+1$ ordinary differential equations
\begin{eqnarray}\label{char}
\dot x= Q(x, {\bY}), \quad \dot \bY = \bS(x, {\bY}).
\end{eqnarray}
One can say that this somewhat simplifies the study of the non-strictly hyperbolic system \eqref{main} due to the fact that one can use results and methods related to the theory of ordinary differential equations.

The interest in studying non-strictly hyperbolic systems is due to the fact that some important multidimensional radially symmetric models of semiconductor and cold plasma physics can be reduced to them. These models are described by the Euler-Poisson equations in the repulsive case with a non-zero constant density background $c>0$,
\begin{equation}\label{1}
\bv_t+ (\bv \cdot \nabla ) \bv =-\nabla \Psi,\quad {\rho}_t+ {\rm div}\,( \rho \bv)
=0,\quad \Delta \Psi =c -\rho,
\end{equation}
where $\bv$, $\rho > 0$, $\Psi$ are the velocity, electron density, and electric field potential, respectively, they depend on time $t\ge 0$ and point $x\in {\mathbb R}^{ d} $.

For
\begin{eqnarray*}\label{sol_form}
\bv=F(t,r) {\br}, \quad \nabla \Psi=\bE=G(t,r) {\br}, \quad
\rho=\rho(t,r),
\end{eqnarray*}
where $ {\br}=(x_1,\dots,x_d)$, $r=|{\br}|$, the system can be reduced to
\begin{eqnarray}\label{sys_pol1}
G_t+F r \,G_ r=F-d F
G, \quad
F_t+F r\, F_r=-F^2 -
G,
\end{eqnarray}
(see \cite{R22_Rad} for details). We see that \eqref{sys_pol1} is a special case of \eqref{main}.
It has a zero constant equilibrium $G=F=0$, and it was natural to suppose that near this equilibrium there are solutions that are also smooth in time.
It is easy enough to show that in dimension one this is indeed the case \cite{ELT}, \cite{RChZAMP21}.

However,  in the multidimensional case the situation is quite different. It has been found that in addition to dimension 1 there is only one spatial dimension 4 for which there is a neighborhood of the zero equilibrium corresponding to globally smooth solutions (the dimensions of this neighborhood were found exactly in \cite{R24_4D}). In the remaining dimensions, any arbitrarily small perturbation of the general form of the equilibrium position leads to a finite-time blow-up.
However, if the perturbation is chosen in a special way so that the initial data lie on some submanifold of lower dimension (corresponding to the so-called {\it simple waves}, such that $F=F(G)$) containing the origin, then a globally smooth solution can still be obtained \cite{R22_Rad}.

As has been recently shown, in all cases where the period of a Lagrangian trajectory $x(t)$ depends on the initial point of the trajectory, the different trajectories necessarily intersect and the solution blows up. Namely, the following lemma holds.
 \begin{lemma}\label{Lem1}
Assume that the mapping $x\mapsto X(t)$ $(\mathbb R \mapsto \mathbb R)$ is continuous, and the trajectory $X(t)$ is (nontrivially) periodic with respect to $t\ge 0$ for all $x_0\in \mathbb R$ with  period $T(x_0)$ that  depends continuously on $x_0$, $X(0)=x_0$. Then, if $T(x)$ is not constant, there exist $x_1$ and $x_2$ from $\mathbb R$ such that  $X_1(t_*)=X_2(t_*)$ for some $t_*>0$, where $X_i(t)$ is the trajectory such that $X_i(0)=x_i$, $i=1,2$.
\end{lemma}
The lemma is proven in \cite{R24_doping}. In fact, this is a more general formulation of Lemma 2.2. \cite{Carrillo}. This result
It is also known to physicists in the context of the problem of destruction of plasma oscillations.

In other words, only in the case where the period of the Lagrangian trajectory $x(t)$ does not depend on the initial point of the trajectory $x(0)$ can one hope to find a neighborhood of zero equilibrium for which globally smooth solutions exist. For problems related to cold plasma physics, cases of isochronous oscillations are extremely rare. Apparently, they are possible only in the non-relativistic case with a constant density background \cite{R24_doping} in dimensions 1 and 4 (for radially symmetric solutions).

However, it is possible to study quasilinear non-strictly hyperbolic systems outside the physical context, that is, to understand what properties $Q$ and $S_i$ must have in order for the \eqref{main} system to have isochronous oscillations.

Let us recall that an ODE system is called {\it isochronous} if its phase space has an open fully dimensional domain, where all its solutions are periodic with the same period, independently from the initial data.

\begin{definition} We will say that the system \eqref{main} is isochronous oscillatory if the characteristic system \eqref{char} is isochronous.
\end{definition}
In particular, this means that each of its characteristics $x(t)$, $x(0)=x_0\in \mathbb R$ is periodic with the same period.

In this paper we first show that an isochronous oscillatory system has a neighborhood of the trivial steady state in the $C^1$ -- norm corresponding to globally time-smooth solutions of the Cauchy problem. We then propose rules that allow generating non-strictly hyperbolic isochronous oscillatory systems. We then discuss the possibility of obtaining smooth solutions for the cold plasma equations in the non-relativistic and relativistic cases.

\section{Global smooth solutions}

We are going to prove the following theorem.
\begin{theorem}\label{T1}
A system \eqref{main} is isochronously oscillatory if and only if there exists a neighborhood $U$ of the point $\bY={\bf 0}$ in the $C^1$--norm such that the solution of the Cauchy problem \eqref{main}, \eqref{mainID} with initial data from $U$ preserves initial smoothness for all $t>0$.
\end{theorem}

\proof  The system \eqref{main} is hyperbolic, so it has a local solution in time, as smooth as the initial data, and the blowup is due either to the unboundedness of the solution components themselves, or to the unboundedness of their first derivatives \cite{Daf16}. The solution components are periodic, so we need to study the behavior of the derivatives. To do this, we differentiate \eqref{main} with respect to $x$ and obtain the following matrix Riccati equation for the vector $\by=(y_1,\dots, y_n)^T=((Y_1)_x,\dots, (Y_n)_x)^T$:
\begin{eqnarray}\label{y}
&&   (y_i)_t + Q \,{\mathbb I} \, (y_i)_x= - (Q_x +\sum\limits_{j=1}^n Q_{Y_j} \, y_j) \,y_i +  \sum\limits_{j=1}^n (S_i)_{Y_j}\, y_i + (S_i)_x,\\&& \quad i,j=1,\dots, n,\nonumber
\end{eqnarray}
with the initial data
$\by^0=(y_1^0,\dots, y_n^0)^T=((Y_1^0)_x,\dots, (Y_n^0)_x)^T$.

The key point in studying the behavior of derivatives is the Radon theorem \cite{Radon}, \cite{Riccati}.

\begin{theorem}\label{TR} [The Radon lemma]
\label{T2} A matrix Riccati equation
\begin{equation}
\label{Ric}
 \dot W =M_{21}(t) +M_{22}(t)  W - W M_{11}(t) - W M_{12}(t) W,
\end{equation}
 {\rm (}$W=W(t)$ is an $(n\times m)$ matrix, $M_{21}$ is an $(n\times m)$ matrix , $M_{22}$ is an $(m\times m)$ matrix, $M_{11}$ is an
 $(n\times n)$ matrix, $M_{12} $ is an $(m\times n)$ matrix{\rm )} is equivalent to the homogeneous linear matrix equation
\begin{equation}
\label{Lin}
 \dot Y =M(t) Y, \quad M=\left(\begin{array}{cc}M_{11}
 & M_{12}\\ M_{21}
 & M_{22}
  \end{array}\right),
\end{equation}
 {\rm (}$Y=Y(t)$  is an $(n\times (n+m))$ matrix, $M$ is an $((n+m)\times (n+m))$  matrix{\rm )} in the following sense.

Let on some interval ${\mathcal J} \in \mathbb R$ the
matrix-function $\,Y(t)=\left(\begin{array}{c}{Q}(t)\\ {P}(t)
  \end{array}\right)$ {\rm (}${Q}$  is an $(n\times n)$ matrix, ${P}$  is an $(n\times m)$ matrix{\rm ) } be a solution of \eqref{Lin}
  with the initial data
  \begin{equation*}\label{LinID}
  Y(0)=\left(\begin{array}{c}{\mathbb I}\\ W_0
  \end{array}\right)
  \end{equation*}
   {\rm (}$ \mathbb I $ is the identity $(n\times n)$ matrix, $W_0$ is a constant $(n\times m)$ matrix{\rm ) } and  $\det {Q}\ne 0$ on ${\mathcal J}$.
  Then
{\bf $ W(t)={P}(t) {Q}^{-1}(t)$} is the solution of \eqref{Ric}   on ${\mathcal J}$   with
$W(0)=W_0$.
\end{theorem}

System \eqref{y} can be re-written along characteristics as \eqref{Ric} for  $W=\by$, with $M_{11}=Q_x $, $M_{12}=(Q_{Y_1},\dots, Q_{Y_n})$, $M_{21}= ((S_1)_x,\dots, (S_n)_x)^T$, and $(n\times n)$ matrix $M_{22}=(S_i)_{Y_j}$, $i,j=1,\dots, n$,   i.e.
\begin{eqnarray}
\label{matr}
 \begin{pmatrix}
 \dot q \\
  \dot y_1\\
  \dots\\
  \dot y_n\\
  \end{pmatrix}
=\begin{pmatrix}
Q_x & Q_{Y_1} &\dots & Q_{Y_n}\\
(S_1)_x& (S_1)_{Y_1}& \dots&(S_1)_{Y_n}\\
\dots&\dots&\dots&\dots\\
(S_n)_x& (S_n)_{Y_1}& \dots&(S_n)_{Y_n}\\
\end{pmatrix}
\begin{pmatrix}
q\\
  y_1\\
  \dots\\
   y_n\\
\end{pmatrix},\quad
\begin{pmatrix}
  q\\
  y_1\\
  \dots\\
   y_n\\
\end{pmatrix}(0)=\begin{pmatrix}
  1\\
   y_1^0\\
  \dots\\
   y_n^0\\
\end{pmatrix}.
\end{eqnarray}

This is a system of linear equations with periodic coefficients that can be studied using Floquet theory (for example, \cite{Chicone}, Section 2.4). According to this theory, for the fundamental matrix $\Psi (t)$ ($\Psi (0)={\mathbb I}$),
 there exists a constant matrix $\mathcal M$, possibly with complex coefficients, such that $\Psi (T)=e^{T \mathcal M}$, where $T$ is the period of the coefficients. The eigenvalues of the matrix of monodromy $e^{T \mathcal M}$ are called the characteristic multipliers of the system. If among the characteristic multipliers there are such that their absolute values are greater than one, then the solution of \eqref{matr} is not periodic and the oscillations of  the components of solution basically have an exponentially rising amplitude (\cite{Chicone}, Theorem 2.53). In particular, this implies that  the zero solution  of \eqref{matr} is unstable in the sense of Lyapunov.
Thus, the component $q$, starting from 1, turns into zero in a finite time, and the solution to \eqref{y} blows up.

If the solution \eqref{matr} is periodic, then
for $y_i^0=0$, $i=1,\dots, n,$ the component $q\equiv 1$ and one can find $(y_1^0, \dots, y_n^0)$ at least close to zero such that $q(t)>0$ within its period. This implies the existence of a neighborhood of zero equilibrium in the $C^1$--norm such that the solution from this neighborhood preserves global smoothness in $t$.

We show that the solution of \eqref{matr} is periodic if and only if \eqref{char} is isochronous.
For this purpose we use the following result, which is in fact a generalization to the multidimensional case of the proposition from \cite{Gasull}, p.8.
\begin{proposition}
Let
\begin{equation}\label{X}
\dot X_i= {\mathcal N}_i({\bf X}), \quad {\bf X}=(X_1,\dots, X_n)^T, \quad i=1,\dots, n, \quad {\mathcal N}=({\mathcal N}_1,\dots, {\mathcal N}_n)^T \in C^2({\mathbb R}),
\end{equation}
be a system with non-negative energy integral $E({\bf X})=h=\rm const$ and there exists a unique point $\overline {\bf X}\in {\mathbb R}^n$ such that $E(\overline {\bf X})=0$. Let ${\bf Z}(t,h)$, ${\bf Z}=(Z_1,\dots, Z_n)^T$ be a periodic solution corresponding to a fixed level $h>0$ with period $T(h)$. Then the solutions of  system
\begin{equation}\label{x}
\dot x_i= \sum\limits_{j=1}^n \,\dfrac{\partial {\mathcal N}_i({\bf Z}(t,h))}{\partial X_j}\, x_j, \quad i=1,\dots, n,
\end{equation}
are $T(h)$ - periodic if and only if $T'(h)=0$.
\end{proposition}
\proof
First of all, we note that ${\bf Z}_h$ solve \eqref{x} with any initial data ${\bf X}^0$ for system \eqref{X}. 
Let us find a condition of periodicity for ${\bf Z}_h$, i.e. ${\bf Z}_h(t,h)={\bf Z}_h(t+T(h),h)$. 

Denote ${\bf Z}'_k$, $k=1,2$, the derivative with respect to the first and second arguments, respectively. Since ${\bf Z}(t,h)={\bf Z}(t+T(h),h)$, then
\begin{equation*}
{\bf Z}_h(t,h)=\dfrac{d \, {\bf Z}(t,h)}{d \, h} =\dfrac{d \, {\bf Z}(t+T(h),h)}{d \, h} ={\bf Z}'_1(t+T(h),h)\, T'(h)+ {\bf Z}'_2(t+T(h),h).
\end{equation*}
 Note that $ {\bf Z}'_1(t+T(h),h)= {\bf Z}_t(t+T(h),h)={\bf \mathcal N}({\bf Z}(t,h))\ne 0$ identically, and   $ {\bf Z}'_2(t+T(h),h)= {\bf Z}_h(t+T(h),h)$. Thus, if  $T'(h)\ne 0$, then ${\bf Z}_h(t+T(h),h)\ne {\bf Z}_h(t,h)$, and vice versa. This proves the proposition. $\Box$

To complete the proof of Theorem \ref{1}, it suffices to note that \eqref{char} can be taken as \eqref{X} with ${\bf X}=(Q, S_1,\dots, S_n)^T$, and \eqref{matr} can be taken as \eqref{x} with ${\bf x}=(q, y_1,\dots, y_n)^T$.
$\Box$

\begin{remark} We emphasize that the requirement of existence of the zero equilibrium position of the system \eqref{main} is essential. Without this requirement, even if the system \eqref{char} has an isochronous equilibrium position, there may be no globally smooth solutions at all. Indeed, consider the Hopf equation with the potential
\begin{equation*}
Y_t+Y Y_x=-x, \quad Y=Y(t,x)
\end{equation*}
The system of characteristics \eqref{char} has the form
\begin{equation*}
\dot x = Y, \quad \dot Y =-x,
\end{equation*}
it coincides with the equations of the harmonic oscillator and has an isochronous center at the origin. Nevertheless, along the characteristics the derivative $y=Y_x$ obeys $\dot y= - y^2-1 $
and, obviously, becomes minus infinity under any initial conditions.
\end{remark}

\section{Construction of isochronous oscillatory systems}

The history of studying isochronous systems of ordinary differential equations is very long. Since the time of Poincar\'e it has been known that if a system has an isochronous equilibrium position, it can be linearized, e.g.\cite{Romanovski}.
However, in practice, finding a transformation that allows a system to be linearized is sometimes much more difficult than finding alternative ways to prove that the system is isochronous. There are  results of this kind for systems on the plane. They can be applied to systems that arise naturally (as happened in the case of system \eqref{sys_pol1}). But if it is difficult to find an isochronous system of natural origin, then if we are talking about designing isochronous systems that have no relation to real physics, then the issue is resolved quite simply.

For example, we can take a system corresponding to a linear oscillator $\dot X_1= X_2, \,\dot X_2= -X_1,$ and consider an invertible smooth transformation $X_1=F_1(Z_1, Z_2), X_2= F_2(Z_1, Z_2)$ such that $F_1(0, 0) = F_2(0,0)=0$. Then the origin of the new system
 \begin{eqnarray*}
  & \dot Z_1= \dfrac{\Delta_1}{\Delta}, \,&\dot Z_2= \dfrac{\Delta_2}{\Delta},\\& \Delta={\rm Det}\begin{pmatrix}
  (F_1)_1& (F_1)_2\\
  (F_2)_1& (F_2)_2\\
  \end{pmatrix},&\quad \Delta_1={\rm Det}\begin{pmatrix}
  F_2& (F_1)_2\\
  -F_1& (F_2)_2\\
  \end{pmatrix},\quad \Delta_2={\rm Det}\begin{pmatrix}
  (F_1)_1& F_2\\
  (F_2)_1& -F_1\\
  \end{pmatrix},
 \end{eqnarray*}
 where $(F_i)_j$ is the derivative of $F_i(Z_1, Z_2)$ with respect to the component $j$, $i,j=1,2$, is the isochronous center.

The same can be done for an uncoupled system of oscillators that are further coupled via an invertible transformation.

Formally, we can choose any isochronous system of ordinary differential equations as \eqref{char} and construct a system \eqref{main} on its basis. If the system constructed in this way has a trivial solution, then it satisfies all the necessary requirements.

For example, the book \cite{Calogero} is devoted to a variety of methods for manufacturing isochronous systems. Thus, if we do not get tied to physics, isochronous systems turn out to be quite numerous, whereas "...in the real world the examples of purely isochronous behavior are rather rare -- otherwise life would be pretty dull" \cite{Calogero}. At the same time, even in the $N$-body problem, it is possible to construct an explicit analytical solution, obtaining isochronous oscillations by special selection of the potential. Among physically meaningful problems let us also note  the cases of isochronicity for  so-called position dependent mass (PDM)-oscillators  \cite{Mustafa}.
 Criteria for the isochronicity of a Hamiltonian system with one degree of freedom was obtained in \cite{Zampieri}, \cite{Treschev}.

\subsection{Homogeneous left-hand sides}

We will focus on a special case, where the system \eqref{char}
is
\begin{eqnarray*}\label{char1}
\dot x= Q(x, {\bY}), \quad \dot \bY = \bS({\bY}).
\end{eqnarray*}
Here the second part $\dot \bY = \bS({\bY})$ is decoupled from the whole system, and the behavior of $x(t)$ is inherited from $\bY$.
This case can be applied to some extent to cold plasma oscillations.

Let the vector $\bY$ consist of two components ($n=2$). In this case we have to study the problem of the center in the plane, and there are many results on isochronicity (see e.g. \cite{Romanovski}, \cite{Kovacic}, \cite{Dong}, \cite{Fernandes}, \cite{Hill} and references therein). In particular, the Sabatini criterion \cite{Sabatini} is very useful for two dimensional systems that can be reduced  in the Lienard equation
\begin{equation}\label{Lienard}
\ddot z+ f (z) \dot z+g(z)=0.
\end{equation}
They are, for example,
\begin{eqnarray*}
 && \dot Y_1=Y_2, \quad \dot Y_2 =-g(Y_1)-f(Y_1) Y_2,
\end{eqnarray*}
and
\begin{eqnarray*}
 && \dot Y_1=Y_2-F(Y_1), \quad \dot Y_2 =-g(Y_1), \quad F(x)=\int\limits_0^x \, f(s)\, ds.
\end{eqnarray*}
The following theorem holds \cite{Sabatini}.
\begin{theorem}\label{S}
Let $f, g$  be analytic, $g$ odd, $f (0)=g(0)=0$, $g'(0)>0$. Then
$\mathcal O =(y,\dot y)=(0,0)$ is a center if and only if $f$ is odd and
$\mathcal O $ is an isochronous center for \eqref{Lienard} if and only if
\begin{equation*}\label{tau}
\tau(z) :=\left(\int\limits_0^z sf(s) ds \right)^2-z^3 (g(z)-g'(0)z)=0.
\end{equation*}
\end{theorem}

\section{Examples}

\subsection{How to stabilize non-relativistic radial oscillations in any dimensions? }

Let us consider system \eqref{sys_pol1}, describing radially symmetric non-relativistic plasma oscillations in $d$-dimensional space.
We rename $r$ to $x$ and $G, F$ to $Y_1, Y_2$, respectively. Thus, system \eqref{char} here is
\begin{eqnarray}\label{sys_pol1_char}
   \dot x = x Y_2,\quad \dot Y_1= Y_2 -d Y_1 Y_2, \quad \dot Y_2= - Y_1 - Y_2^2,
     \end{eqnarray}
the latter two equations are reduced to
\begin{equation*}
\ddot Y_2+(2+d)\,Y_2\,\dot  Y_2 + Y_2+d \, Y_2^3=0,
\end{equation*}
which is the Lienard equation, so we can apply Theorem \ref{S}. We see that
\begin{equation*}\label{tauY}
\tau(
Y_2) := ((2+d)^2- 9 d) \,\frac{Y_2^6}{2},
\end{equation*}
therefore $\tau(Y_2)=0 $ implies $d=1$ or $d=4$.

 We study the following question: what velocity-dependent force ${\bf F} (\bV, r)$ should be added to the right-hand side of the first equation \eqref{1} to ensure that the oscillations are isochronous and hence that globally smooth perturbations of the trivial steady state can exist in any dimension? It is easy to see that we need to calibrate the coefficients of the system \eqref{sys_pol1_char} depending on $d$.

First of all, we denote $L(x,Y_2)=\frac{1}{r} {\bf F} (\bV, r)|_{r=x}=\frac{1}{r} {\bf F} (r Y_2, r)|_{r=x}$. Then \eqref{sys_pol1_char}
takes the form
\begin{eqnarray*}\label{sys_pol1_charM}
   \dot x = x Y_2,\quad \dot Y_1= Y_2 -d Y_1 Y_2, \quad \dot Y_2= - Y_1 - Y_2^2-L(x,Y_2).
   \end{eqnarray*}
   Standard computations show that to obtain the Lienard equation we have to require $L=L(Y_2)$ and if we hope to obtain $\tau(Y_2)=0$,
   we have to set $L=\gamma Y_2^2$, $\gamma=\rm const$. In this case,  $\tau(Y_2)=0 $, if $\gamma$ is a root of the quadratic equation
\begin{equation*}
(2(1-\gamma)+d)^2- 9 d(1-\gamma)=0,
\end{equation*}
that is, $\gamma=1-d$ or $\gamma=1-\frac{d}{4}$. The respective force term is
\begin{equation*}
 {\bf F} (\bV, r)=
\gamma \,\dfrac{|\bV|^2}{r}.
\end{equation*}

It is tempting to interpret this term as something like aerodynamic friction, but this is not the case, since friction is directed against the velocity and must be proportional to $ - \bV |\bV|$. Moreover, the presence of friction is always associated with a decay in the energy integral, which is incompatible with the existence of a center on a plane $(G,F)$. In reality, the calibration term has no physical meaning and only means that the velocity must be specially slowed down or increased at different stages of the oscillation.

\subsection{Is it possible  to stabilize relativistic oscillations in $1D$?}

The system of Euler-Poisson equations describing the behavior of relativistic cold plasma in $ {\mathbb R}$  with a variable density background in the repulsive case has the following form \cite{book}:
\begin{equation*}\label{1R}
P_t+ (V \cdot \nabla ) P =-E,\quad \dfrac{\partial \rho}{\partial t}+ {\rm div}\,(\rho V)
=0,\quad E =c(x) -\rho,\quad  V = \dfrac{P}{\sqrt{1+P^2}}.
\end{equation*}
Here $x$ and $t$ are dimensionless coordinates in
space and time, respectively. The variable $P$ describes
the electron momentum, $V$ is the electron velocity, $E$ is a function
characterizing the electric field,
$\rho> 0$ is electron density.   A fixed $C^1$ - smooth function
$c(x)>0$ is the density background or the so-called doping profile. In the simplest case $c(x)=1$.

Using it, we arrive at equations describing
flat one-dimensional relativistic plasma oscillations 
\begin{equation}\label{u1}
P_t +
V\,P_x  + E = 0,\quad
E_t  +
V\, E_x - V c(x)= 0,
\quad V = \dfrac{P}{\sqrt{1+P^2}}.
\end{equation}
 with initial conditions
\begin{equation}\label{cd1}
P(x,0) = P_0(x), \quad
E(x,0) = E_0(x), \quad
x \in {\mathbb R}.
\end{equation}

Thus, along characteristics $x=x(t)$, starting from a point $x_0\in \mathbb R$ the solution $(V(x(t)),E(x(t)))$ obeys the system of ODEs
\begin{equation}\label{char_solR}
\dot x=\frac{P}{\sqrt{1+P^2}}, \quad \dot P =-E, \quad \dot E  =c(x) \frac{P}{\sqrt{1+P^2}},
\end{equation}
with the initial data
\begin{equation*}\label{char_sol_CD}
x(0)=x_0, \quad P(0)=P_0(x_0), \quad E(0)=E_0(x_0).
\end{equation*}

We prove several proposition allowing to understand properties of plain relativistic oscillations.
\begin{proposition}\label{P1}
In the case of a constant doping profile $c(x)=c_0>0$ any nontrivial classical solution to the Cauchy problem \eqref{u1}, \eqref{cd1}, which is not a simple wave $P=P(E)$, blows up in a finite time.
\end{proposition}

\proof  From two latter equations of
 \eqref{char_solR} we have
 \begin{equation} \label{LP1}
   \ddot P + c_0 \frac{P}{\sqrt{1+P^2}}=0,
 \end{equation}
which is a particular case of the Li\'enard equation \eqref{Lienard}. It is easy to see that $\tau(P)\ne 0$, therefore Theorem \ref{S}
implies that the oscillations of $P$ are not isochronous. From the first equations \eqref{char_solR} we conclude that the oscillations of $x(t)$ are also not isochronous. Thus, from Lemma \ref{Lem1} we conclude that the characteristics necessarily intersect and the solution blows up on a finite time for general initial data.
In \cite{RChZAMP21} it is shown that for a simple wave one can choose a neighborhood of the origin $P=G=0$ corresponding to a smooth solution. $\Box$

\medskip

We consider an analog of system \eqref{1} with a forcing, such that the first equation takes the form
\begin{equation*}
 P_t+ (V \cdot \nabla ) P =-E + {\mathcal L},
\end{equation*}
where ${\mathcal L}$ is an exterior force. We assume  ${\mathcal L}'(0)=0$, ${\mathcal L}'(P)$ is odd, to preserve the oscillating character of the forced system.

\begin{proposition} No momentum (and velocity) dependent force  ${\mathcal L}(P)$ such that ${\mathcal L}'(0)=0$, ${\mathcal L}'(P)$ is odd, can help to make the oscillations of \eqref{1} isochronous in the case of a constant doping profile.
\end{proposition}

\proof Let us consider the respective analog of \eqref{LP1},
 \begin{equation*} 
   \ddot P + {\mathcal L}'(P) \dot P +  c_0 \frac{P}{\sqrt{1+P^2}}=0.
 \end{equation*}

 Thus,
 \begin{equation*}
   \tau(P)=\left(\int\limits_0^P \, s {\mathcal L}'(s)\,ds\right)^2 - c_0 \, P^4 \,\frac{1-\sqrt{1+P^2}}{\sqrt{1+P^2}}.
 \end{equation*}
Since $1-\sqrt{1+P^2}<0$ and $c_0>0$, then $\tau>0$ for all nontrivial $P$. $\Box$

\bigskip

In \cite{R24_doping} it was proved that in the case of non-relativistic cold plasma oscillations the system os characteristics is isochronous (and
therefore globally smooth solutions are possible only for a constant doping profile. Further, as follows from Proposition \ref{P1}, relativistic  cold plasma oscillations for a constant doping profile basically blow up. Therefore, the question arises: is it possible to find a variable $c(x)$ such that the oscillations do not blow up? In other words, can the system \eqref{char_solR} be isochronous for some choice of $c(x)$? The next proposition states that the answer is negative.

\begin{proposition}
 For any smooth doping profile $c(x)>0$, $x\in\mathbb R$, any nontrivial classical solution of the Cauchy problem \eqref{u1}, \eqref{cd1} blows up in a finite time.
\end{proposition}
\proof
From \eqref{char_solR} we have
\begin{eqnarray*}
E(x(t)) =\int\limits_{x_0}^x \, c(\xi)\, d\xi + E_0(x_0), \quad P^2= \frac{(\dot x)^2}{1-(\dot x)^2}, \quad \ddot x=\frac{1}{(1+P^2)^\frac{3}{2}}\,\dot P,
\end{eqnarray*}
therefore
\begin{equation}\label{ddotx1}
\ddot x=-E(x)\,(1-(\dot x)^2)^\frac{3}{2}.
\end{equation}
Let us notice that $|\dot x|< 1$. Assume for the sake of simplicity that the equilibrium point is $x=0$.

If we denote $\dot x = s(x)$ and $z(x)=s^2$, we get
\begin{equation*}
z'=-2 E(x)\,(1-z)^\frac{3}{2},
\end{equation*}
and
\begin{eqnarray*}
&&1-(\dot x)^2=1-z=\frac{1}{\Phi^2(x)}, \quad \Phi(x)=-\int\limits_0^x E(\xi) d\xi +E_0+(1+P_0^2)^\frac{3}{2},\\ && \Phi'(x)=-E(x), \quad \Phi''(x)=c(x).
\end{eqnarray*}
Together with \eqref{ddotx1} it implies
\begin{equation}\label{ddotx}
\ddot x=\frac{\Phi'(x)}{\Phi^3(x)},
\end{equation}
which corresponds to a plane Hamiltonian system with Hamiltonian $\mathcal V$ such that $g(x)={\mathcal V}'=-\frac{\Phi'(x)}{\Phi^3(x)}$.

1. If we assume that $g(x)$ is {\it odd} and analytic in the neighborhood of $x=0$, then we can use Theorem \ref{S} and calculate
\begin{equation*}
\tau(x)=-x^3 \left(-\frac{\Phi'(x)}{\Phi^3(x)}-K\, x\right), \quad K= -\,\left(\frac{\Phi'(x)}{\Phi^3(x)}\right)'\Big|_{x=0}>0.
\end{equation*}
Thus, $\tau(x)=0$ if and only if
\begin{equation*}\label{Phi}
\Phi(x)=\,\pm \,\frac{1}{\sqrt{K \, x^2 +M}}, \quad M={\rm const}.
\end{equation*}
Note that in this case $\frac{\Phi'(x)}{\Phi^3(x)}=-K x$, therefore  \eqref{ddotx} is linear and it is invariant invariant under the shift $x_0$.
 Considering the shift of the initial point to $x_0$, we have
\begin{equation*}
 c(x)=\pm \, \frac{K (M-2K (x-x_0)^2)}{(M+K (x-x_0)^2)^{\frac{5}{2}}}.
\end{equation*}

 However, we see that $c(x)$ cannot be positive for all $x\in \mathbb R$. Moreover, to ensure isochronicity $c(x)$ must depend on the initial point of the trajectory, which contradicts the requirement that the doping profile be a function of the spatial coordinate only.

 2. It is known that the assumption that $g(x)$ is odd is not necessary for the existence of an isochronous center of the Hamiltonian system. Thus, we can use a general result of \cite{Zampieri} for the class of continuous functions $g$.

 A $C^1$ - diffeomorphism $H$ of an open interval $J  \subseteq {\mathbb R}$
onto itself is called an {\it involution} if
\begin{equation*}
  H^{-1} = H,\quad  0 \in J,\quad  H(0) = 0, \quad H'(0) = -1.
\end{equation*}
It means that the graph
of $y=H(x)$ is symmetric with respect to the main diagonal $y=x$.

\begin{theorem}\label{Z}
 Let $H : J \rightarrow J$ be an involution, $\omega > 0$, and define
\begin{equation}\label{2.4}
  V (x) =
\frac{\omega}{8} \,(x - H(x))^2,\quad x \in J.
\end{equation}

Then the origin is an isochronous center for $\ddot x = -g(x),$ where $g(x) =
V '(x),$ with the same period $\dfrac{2\pi}{\omega}$. Vice
versa, let $g$ be continuous on a neighbourhood of $0 \in {\mathbb R},$ $ g(0) = 0,$
suppose there exists $g'(0) > 0$, and the origin is an isochronous center
for $\ddot x = -g(x)$, then there exist an open interval $J,\,  0 \in J$, which is a
subset of the domain of $g$, and an involution $H : J \rightarrow J$ such that \eqref{2.4}
holds with
\begin{equation*}
 V (x) = \int\limits_0^x g(s)ds, \qquad  \omega =
\sqrt{g'(0)}.
\end{equation*}
\end{theorem}

The paper \cite{Zampieri} contains a lot of examples for $H$. In our case \eqref{ddotx} we have
\begin{eqnarray*}
&&  M-\frac{1}{2 \Phi^2(x)}= -\,\frac{\omega}{8} \,(x - H(x))^2, \\ && M= \frac{1}{2 \Phi^2(0)}>0,\quad \Phi(x)=\pm \, \frac{1}{\sqrt{2M+ \frac{\omega}{4} \,(x- x_0 - H(x-x_0))^2}}.
\end{eqnarray*}
Thus, in this case $c(x)=\Phi''(x)$ again depends on the starting point of the trajectory and does not satisfy the requirements. Thus the proof is complete.
 $\Box$

\section{Discussion}

1. If the oscillatory system is isochronous, this actually means that it has an additional first integral. In the context of Hamiltonian systems, this property is called superintegrability and has numerous applications \cite{Miller}.

2. The methods presented here can be extended to the systems of partial differential equations for $\bY=\bY(t, \bx)$, $\bx \in {\mathbb R}^m$ of the form
\begin{eqnarray}\label{new}
 & (Y_i)_t+\sum\limits_{j=0}^m a_j(\bY, \bx) (Y_i)_{x_j} =S_i (\bY, \bx), \\\nonumber &\bx=(x_1,\dots, x_m),\quad 
  {\bf a}=(a_1,\dots, a_m),\quad {\bf a}({\bf 0}, \bx)={\bf 0}, \quad
 i=1,\dots, n,\,j=1,\dots, m.
\end{eqnarray}
In particular, for $n=1$, it is one equation. The dynamics along characteristics is defined by the system of $n+m$ ordinary differential equations
\begin{eqnarray}\label{charM}
\dot \bx= {\bf a}(\bx, {\bY}), \quad \dot \bY = \bS(\bx, {\bY}).
\end{eqnarray}
If the equilibrium of \eqref{charM} is isochronous, then it is likely that globally smooth solutions \eqref{new} can be found near the zero steady state.

\section*{Acknowledgements}

The author is grateful to V.V.Bykov for discussing different aspects of isochronous oscillations.


\end{document}